\documentclass[12pt]{article}

\usepackage{amssymb}
\usepackage{amsxtra}
\usepackage{amsmath}
\usepackage{amsfonts}
\usepackage{CJK}
\usepackage[titletoc]{appendix}
\usepackage{graphicx}
\usepackage{color}
\numberwithin{equation}{section}

\newtheorem{thm}{Theorem}[section]

\newtheorem{lem}[thm]{Lemma}

\newtheorem{rem}[thm]{Remark}

\newcommand{\eqa}{\begin{eqnarray}}
\newcommand{\eeqa}{\end{eqnarray}}
\newcommand{\beq}{\begin{equation}}
\newcommand{\eeq}{\end{equation}}
\newcommand{\nn}{\nonumber}
\newcommand{\p}{\partial}

\def \dsum{\displaystyle\sum}

\newcommand{\diag}{\mathrm{diag}}

\allowdisplaybreaks
\begin{document}
\title{ {The Novel Symmetry Constraint and Binary Nonlinearization of the Super Generalized Broer-Kaup Hierarchy with Self-consistent Sources and Conservation Laws
\thanks{\footnotesize The work was partially supported by the National Natural Science Foundation of
China under Grant No. 11601055, Natural Science Foundation of Anhui Province under
Grant No. 1408085QA06.}}}

\author{{\footnotesize { Bei-bei Hu$^{1}$ ,Fang Fang$^{1,}$\thanks{Corresponding authors. E-mails: hubsquare@chzu.edu.cn; fangfang7679@163.com; zhangningsdust@126.com} ,  Ning Zhang$^{2}$ }}\\
{\footnotesize { {\it$^{1}$ School of Mathematics and Finance, Chuzhou University, Anhui 239000, China}}}\\
{\footnotesize { \it $^{2}$ Department of Basical Courses, Shandong University of Science and Technology, Taian 271019, China}}
}
\date{\small }\maketitle

\textbf{Abstract}: The super generalized Broer-Kaup(gBK) hierarchy and its super Hamiltonian structure
are established based on a loop super Lie algebra and super-trace identity.
Then the self-consistent sources, the conservation laws,
the novel symmetry constraint and the binary nonlinearization
of the super gBK hierarchy are generated, respectively.
In addition, the integrals of motion required for Liouville integrability are
explicitly given.\\
\\
\textbf{ PACS numbers} {02.30.Ik, 02.30.Jr, 02.20.Sv}\\
\textbf{Keywords}: {Super generalized Broer-Kaup hierarchy, Super Hamiltonian structure, Self-Consistent Sources, Conservation Laws, Binary nonlinearization}

\section{ Introduction}

It is known that super integrable systems provide interesting and important models in the supersymmetry theory
Supersymmetry is originated in 1970s when physicists have proposed simple models with
supersymmetric colors in string models and mathematical physics respectively. After that,
Wess and Zumino \cite{wess1974} applied supersymmetry to the four-dimensional spacetime. Unfortunately,
the supersymmetry partners of any particle have not been found so far, and it is generally
believed that this symmetry is spontaneous rupture. In order to unify two kinds of particles
with different spin and statistical properties-Boson and Fermion, theoretical
physicists proposed the concept of hyperspace in the study of unified field theory and quantum
field theory. Inspired by this, mathematicians developed the super analysis, the hypergeometric and the super algebra.

Due to the importance of supersymmetry in physics(especially in the exploration of the relationship
between supersymmetric conformal field and chord theory), which has captured great attention
for the work of super integrable systems associated with Lie super algebra \cite{Kupershmidt1987}, a multitude of classical integrable equations
have been extended to be the super completely integrable equations (see Refs [3-11]
and references therein). Among those, Hu \cite{hu1990,hu1997} and Ma \cite{ma2008} has made a great work. In 1990, Hu proposed the
super-trace identity in his Ph.D. thesis \cite{hu1990}, which is an effective tool to constructing super
Hamiltonian structures of super integrable systems. In 2008, Ma given a systematic proof
of the super-trace identity and the super double Hamiltonian structure
of many super integrable equations is established by using of the super-trace identity
(see Refs \cite{ma2008} and references therein).

Soliton equation with self-consistent sources is an important part in soliton theory.
They are relevant to some problems related to hydrodynamics,
solid state physics, plasma physics, and they are also usually used to describe
interactions between different solitary waves, such as the NLS equation with self-consistent
sources can describe the propagation of solitary waves in the medium of resonance and non resonant media.
It can also describe the interaction between high frequency static wave and ion acoustic wave in plasma \cite{Doktrov1983},
the KP equation with self-consistent sources description the interaction of between short wave and length wave spread
in the X-Y plane \cite{Zakharov1986}, the KdV equation with self-consistent sources description of the interaction of the interaction
of plasma high-frequency wave packet weight and a low frequency wave packets \cite{Mernikov1990}.
And the conservation laws is also an important part in soliton theory.
An infinite number of conservation laws for KdV equation were first discovered by Miura\cite{Miura1968} et al. in 1968,
and then lots of methods have been developed to find them. This may be mainly due to the contributions of
Wadati and others [19-21]. Conservation laws also play an important role in mathematics and engineering as well.
Many papers dealing with symmetries and conservation laws were presented. The direct construction method of
multipliers for the conservation laws was presented \cite{Bluman2002}.
Thus, the study of integrable equations with self-consistent sources and conservation laws has received
much attention. Recently, an army of classical integrable equations to be the super integrable equations,
by using symmetry constraints, and the super integrable system with self-consistent sources and
conservation laws of the super integrable system are constructed(see Refs[22-30] and references therein).

In recent 10 years, He et al. applied the binary nonlinearization method to the super integrable systems (see Refs [30-35]
and references therein) to the construction of finite dimensional super integrable systems from super
soliton equations by using symmetry constraints. It is well known that a crucial idea in carrying out
symmetry constraints is the nonlinearization of Lax pairs for soliton hierarchies. The nonlinearization
of Lax pairs can be classified into mono-nonlinearization, which is proposed by Cao\cite{cao1992} and
binary nonlinearization, which is proposed by Ma [21,36-39] and has attracted a lot of interest recently\cite{zhou2011}.

In 2013, Zhang, Han and Tam \cite{Zhang2013} making use of a Lie algebra and Tu-Ma scheme obtained a new generalized
Broer-Kaup (gBK) equation, where the spatial spectral problem is given by
\eqa {\small \phi_x=M\phi,\quad M=\left(\begin{array}{cc}
-\lambda+\frac{v}{2}&1\\
-2w-2&\lambda-\frac{v}{2}
\end{array} \right)},\label{1.1}\eeqa
where $v$ and $w$ are both scalar potentials, $\lambda$ is the spectral parameter,
and in \cite{Zhang2013}, they have presented two kinds of Darboux transformations, the bilinear presentation,
the bilinear B\"{a}cklund transformation and the new Lax pair of the gBK equation, respectively, by employing the Bell polynomials.
In this paper, we consider the conservation laws, self-consistent sources and the binary nonlinearization
of the super gBK hierarchy under the novel symmetry constraint.
In addition, under the symmetry constraint, the $n$-th flow for a super gBK hierarchy is decomposed
into two super finite-dimensional integrable Hamilton systems over the super symmetric manifold.

Organization of this paper. In the next section, we construct the super gBK hierarchy and its super
Hamiltonian structure based on a loop super Lie algebra and super-trace identity.
In Section 3, we construct the super gBK hierarchy with self-consistent sources.
And conservation laws of the super gBK hierarchy are constructed in Section 4.
In Section 5, we obtain a symmetry constraint for the potential of the super gBK hierarchy.
Then in Section 6, we apply the binary nonlinearization to the super gBK equation hierarchy using a symmetry
constraint for the potential of the super gBK hierarchy and obtain a super finite-dimensional
integrable Hamiltonian system on the super symmetry manifold, whose integrals of motion are explicitly
given. And some conclusions are given in the last Section.

\section{ The super gBK hierarchy }

In this section, we shall construct a super gBK hierarchy starting from a Lie super-algebra.
We consider the following spatial spectral problem
 \eqa {\small \phi_x=M\phi,\quad M=\left(\begin{array}{ccc}
-\lambda+\frac{v}{2}&1&\alpha\\
-2w-2&\lambda-\frac{v}{2}&\beta\\
\beta&-\alpha&0
\end{array} \right),\quad
\phi=\left(\begin{array}{c}
\phi_1\\
\phi_2\\
\phi_3\end{array} \right),\quad u=\left(\begin{array}{c}
v\\
w\\
\alpha\\
\beta\end{array} \right)},\label{2.1}\eeqa
where $\lambda$ is the spectral parameter, $v$ and $w$ are even potentials, and $\alpha$ and $\beta$ are odd potentials.

And associated with the Lie superalgebra $sl(2,1)$. Its basis is
\eqa &&
e_1=\left(\begin{array}{ccc}
1&0&0\\
0&-1&0\\
0&0&0\end{array} \right),
e_2=\left(\begin{array}{ccc}
0&0&0\\
1&0&0\\
0&0&0\end{array} \right),
e_3=\left(\begin{array}{ccc}
0&1&0\\
0&0&0\\
0&0&0\end{array} \right),\nn\\
&&e_4=\left(\begin{array}{ccc}
0&0&1\\
0&0&0\\
0&-1&0\end{array} \right),
e_5=\left(\begin{array}{ccc}
0&0&0\\
0&0&1\\
1&0&0\end{array} \right).\nn \eeqa
where $e_1,e_2,e_3$ are even elements and $e_4,e_5$ are odd ones,
$[.,.]$ and $[.,.]_{+}$ denote the commutator and the anticommutator,
satisfy the following operational relations:
\eqa &&[e_1,e_2]=-2e_2,[e_1,e_3]=2e_3,[e_2,e_3]=-e_1,\nn\\
&&[e_5,e_1]=[e_2,e_4]=e_5,[e_3,e_4]=[e_2,e_5]=0,[e_3,e_5]=[e_1,e_4]=e_4,\nn\\
&&[e_4,e_4]_{+}=-2e_3,[e_5,e_5]_{+}=2e_2,[e_4,e_5]_{+}=[e_5,e_4]_{+}=e_1. \label{2.2}\eeqa
Then corresponding loop super algebra is given by
\beq sl(2,1)=sl(2,1)\otimes \mathbb{C}[\lambda,\lambda^{-1}],\nn\eeq
where $\mathbb{C}[\lambda,\lambda^{-1}]$ presents the set of Laurent polynomials in $\lambda$ over the complex
number set $\mathbb{C}$.

The corresponding (anti)commutative relations are given as
\beq [e_i\lambda^m,e_j\lambda^n]=[e_i,e_j]\lambda^{m+n},\forall e_i,e_j\in sl(2,1).\nn\eeq

From the Tu format, we setting
\eqa {\small N=\left(\begin{array}{ccc}
A&B&\rho\\
-2B+2C&-A&\delta\\
\delta&-\rho&0\end{array} \right)=\sum_{m\geq
0}\left(\begin{array}{ccc}
a_m&b_m&\rho_m\\
-2b_m+2c_m&-a_m&\delta_m\\
\delta_m&-\rho_m&0\end{array} \right)\lambda^{-m} },\label{2.3}\eeqa
the corresponding $A,B,C$ are even elements and $\rho,\delta$ are odd elements, if we want to get the super integrable system,
we solve the stationary zero curvature equation at first
 \beq N_x=[M,N].\label{2.4}\eeq
Substituting $M, N$ into Eq.\eqref{2.4} and comparing the coefficients of $\lambda^{-m} (m \geq 0)$, we have
\eqa
\left\{\begin{array}{l}
(a_{m+1}, -2b_{m+1},-2\delta_{m+1}, 2\rho_{m+1})^T=\mathcal{L}( a_{m} ,-2b_{m},-2\delta_{m}, 2\rho_{m})^T,\\
a_m=\p^{-1}(2wb_m+2c_m+\alpha\delta_m+\beta\rho_m).\\
\end{array}\right.\label{2.5}\eeqa
Where the recursion operator $\mathcal{L}$ has the following form
\beq \mathcal{L}=\left(\begin{array}{cccc}
\frac{1}{2}\p^{-1}v\p+\frac{1}{2}\p& \frac{1}{2}\p^{-1}w\p+\frac{1}{2}w+1&
-\frac{1}{2}\p^{-1}\alpha\p+\frac{1}{4}\alpha& -\frac{1}{2}\p^{-1}\beta\p-\frac{1}{4}\beta\\
2& -\frac{1}{2}\p+\frac{1}{2}v  & 0 &\alpha\\
2\beta-2\alpha\p & -2\alpha(w+1) &\p+\frac{1}{2}v &\alpha\beta-(2w+2) \\
-2\alpha &\beta &-1 & -\p+\frac{1}{2}v \end{array}
\right).\label{2.6}\eeq

For a given initial value $a_0=k_0\neq0,b_0=c_0=\rho_0=\delta_0=0$, the $a_j,b_j,c_j,\rho_j,\delta_j$$(j\geq1)$ can be
calculated by the recursion relation \eqref{2.5}. Here we list the several values
\eqa &&a_1=0,b_1=-k_0,c_1=k_0w,\rho_1=-k_0\alpha,\delta_1=-k_0\beta,a_2=k_0w-k_0\alpha\beta,\nn\\
&&b_2=-\frac{1}{2}k_0v,c_2=\frac{1}{2}k_0w_{x}+\frac{1}{2}k_0vw,\rho_2=k_0\alpha_x-\frac{1}{2}k_0v\alpha,\delta_2=-k_0\beta_x-\frac{1}{2}k_0v\beta,\nn\\
&&a_3=k_0(\frac{1}{2}w_x+v+wv-v\alpha\beta+\alpha_x\beta-\alpha\beta_x),b_3=k_0(\frac{1}{4}v_x-w-\frac{1}{4}v^2-\alpha\alpha_x+\alpha\beta),\nn\\
&&c_3=\frac{1}{4}k_0(w_{xx}+(wv)_x+vw_x+wv^2)+\frac{1}{2}k_0(v_x+2w^2)-k_0(w\alpha\beta+\alpha\alpha_x-\frac{1}{2}\beta\beta_x),\nn\\
&&\rho_3=k_0(-\alpha_{xx}+v\alpha_x+\frac{1}{2}\alpha v_x-\frac{1}{4}\alpha v^2-\alpha w-\beta_x),\nn\\
&&\delta_3=k_0(-\beta_{xx}-v\beta_x-\frac{1}{2}\beta v_x-\frac{1}{4}\beta v^2-\beta w+(2w+2)\alpha_x+w_x\alpha). \nn\eeqa
Then, consider the auxiliary spectral problem associated with the spectral problem \eqref{2.1}
$$\phi_{t_n}=N^{(n)}\phi$$
where
\beq N^{(n)}=N_{+}^{(n)}+\Delta_n=\sum^n_{m=0}\left(\begin{array}{ccc}
a_m&b_m&\rho_m\\
-2b_m+2c_m&-a_m&\delta_m\\
\delta_m&-\rho_m&0\end{array}
\right)\lambda^{n-m}+\left(\begin{array}{ccc}
b_{n+1}&0&0\\
0&-b_{n+1}&0\\
0&0&0\end{array}
\right)
. \label{2.7} \eeq
with $\Delta_n$ being the modification term, substituting Eq.\eqref{2.7} into the zero curvature equation
\beq U_t-N_x^{(n)}+[U,N^{(n)}]=0, \label{2.8} \eeq
we can obtain the following super gBK hierarchy
\beq
{\small u_{t_n}=\left(\begin{array}{c}
v\\
w\\
\alpha\\
\beta\end{array} \right)_{t_n}=\left(\begin{array}{c}
2b_{n+1,x}\\
-a_{n+1,x}+\alpha\delta_{n+1}+\beta\rho_{n+1}\\
\alpha b_{n+1}-\rho_{n+1}\\
-\beta b_{n+1}+\delta_{n+1}\end{array} \right)=J\left(\begin{array}{c}
a_{n+1}\\
-2b_{n+1}\\
-2\delta_{n+1}\\
2\rho_{n+1}\end{array} \right)}. \label{2.9} \eeq
where the super-Hamiltonian operator $J$ is given by
\beq
J=\left(\begin{array}{cccc}
0 &-\p &0&0\\
-\p &0 &-\frac{1}{2}\alpha&\frac{1}{2}\beta\\
0 &-\frac{1}{2}\alpha &0&-\frac{1}{2} \\
0 &\frac{1}{2}\beta &-\frac{1}{2}&0 \end{array}
\right)\label{2.10}. \eeq
Taking $k_0=2,n=2,t_2=t$, the Eq.\eqref{2.9} can be reduced to the super gBK equation
 \beq  \left\{
\begin{array}{l}
v_{t}=v_{xx}-2vv_x-4w_x-4\alpha\alpha_{xx}+4\alpha_x\beta+4\alpha\beta_x,\\
w_{t}=-w_{xx}-2(wv)_x-2v_x+2(2w+2)\alpha_x+2w_x\alpha-(2w\beta+\frac{1}{2}\beta v^2)(1+\alpha)-2\beta\beta_x,\\
\alpha_{t}=2\alpha_{xx}-2v\alpha_x+2\beta_x-\frac{1}{2}\alpha v_x,\\
\beta_{t}=-2\beta_{xx}-2v\beta_x+2\beta\alpha\alpha_x+2(2w+2)\alpha_x+2\alpha w_x-\frac{3}{2}\beta v_x,
\end{array}
\right. \label{2.11} \eeq
whose Lax pair are $M$ and $N^{(2)}$, $N^{(2)}$ has the following form
\beq
N^{(2)}=\left(\begin{array}{ccc}
2\lambda^2+\frac{1}{2}(v_{x}-v^2)-2\alpha\alpha_x & -2\lambda-v & 2\alpha\lambda+2\alpha_x-v\alpha\\
4(1+w)\lambda+2w_x+2v(1+w) & -2\lambda^2-\frac{1}{2}(v_{x}-v^2)+2\alpha\alpha_x & -2\beta\lambda-2\beta_x-v\beta\\
-2\beta\lambda-2\beta_x-v\beta & -2\alpha\lambda-2\alpha_x+v\alpha &0
\end{array} \right).\label{2.12}\eeq
when $\beta=\alpha=0$, Eq.\eqref{2.8} just reduces to the gBK equation
\eqa \left\{ \begin{array}{l}
v_{t}=v_{xx}-2vv_x-4w_x,\\
w_{t}=-w_{xx}-2(wv)_x-2v_x.
\end{array}
\right. \label{2.13} \eeqa

Next, we use the super trace identity, which proposed by Hu in \cite{hu1997}
and rigorously proved by Ma et al. in ref.\cite{ma2008}:
\beq \frac{\delta}{\delta u}\int Str(N\frac{\p M}{\p \lambda})dx=(\lambda^{-\gamma}\frac{\p}{\p\lambda}\lambda^{\gamma})Str(N\frac{\p M}{\p u}),\label{2.14}\eeq
where $Str$ denotes the super trace. It is not difficult to find that
\eqa &&
Str(N\frac{\p M}{\p \lambda})=-2 A,
Str(N\frac{\p M}{\p v})=A,
Str(N\frac{\p M}{\p w})=-2B,\nn\\&&
Str(N\frac{\p M}{\p \alpha})=-2\delta,
Str(N\frac{\p M}{\p \beta})=2\rho,
\label{2.15}\eeqa
substituting Eq.\eqref{2.15} into Eq.\eqref{2.14}, and comparing the coefficient of $\lambda^{-n-1}$ of both sides of Eq.\eqref{2.14}, we have
\beq \left(\begin{array}{c}
\frac{\delta}{\delta v}\\
\frac{\delta}{\delta w}\\
\frac{\delta}{\delta \alpha}\\
\frac{\delta}{\delta \beta}\end{array}
\right)\int -2a_{n+1}dx=(\gamma-n)\left(\begin{array}{c}
a_{n}\\
-2b_{n}\\
-2\delta_{n}\\
2\rho_{n}\end{array} \right). \label{2.16}\eeq
To fix the constant $\gamma$, we set $n=1$ in \eqref{2.16} and find that
$\gamma=0$, Thus we have
\beq
\frac{\delta \tilde{H}_n}{\delta u}=\left(\begin{array}{c}
a_{n}\\
-2b_{n}\\
-2\delta_{n}\\
2\rho_{n}\end{array} \right),
\tilde{H}_n=\int\frac{2a_{n+1}}{n}dx, n\geq0,
\label{2.17}\eeq
specially, by making use of the recursive relationship \eqref{2.5}, the super gBK equation hierarchy \eqref{2.9} possesse the
following super bi-Hamiltonian structure
\beq u_{t_n}=J\left(\begin{array}{c}
a_{n+1}\\
-2b_{n+1}\\
-2\delta_{n+1}\\
2\rho_{n+1}\end{array} \right)=J\mathcal{L}\left(\begin{array}{c}
a_{n}\\
-2b_{n}\\
-2\delta_{n}\\
2\rho_{n}\end{array} \right)=J \mathcal{L}\frac{\delta \tilde{H}_n}{\delta u}, n\geq0,\label{2.18}\eeq
where the second super-Hamiltonian operator $R$ is given by
\beq
R=J\mathcal{L}=\left(\begin{array}{cccc}
-2\p & \frac{1}{2}\p^2-\frac{1}{2}\p v & 0 & \p\alpha \\
-\frac{1}{2}v\p-\frac{1}{2}\p^2 & -w\p-\frac{1}{2}\p w-\p & R_{1}
& R_{2} \\
0 & \frac{1}{4}\alpha\p-\frac{1}{4}v\alpha-\frac{1}{2}\beta  & \frac{1}{2}  &\frac{1}{2}\p-\frac{1}{4}v  \\
\alpha\p & -\frac{1}{4}\beta\p+(w+1)\alpha+\frac{1}{4}v\beta & -\frac{1}{2}\p-\frac{1}{4}v & (w+1)-\alpha\beta \end{array}
\right)\label{2.19} \eeq
with
\eqa &&
R_{1}=-\frac{1}{4}\p\alpha-\frac{1}{4}v\alpha-\frac{1}{2}\beta,\quad
R_{2}=\frac{1}{4}\p\beta+(w+1)\alpha+\frac{1}{4}v\beta.\nn\eeqa

\section{The super gBK hierarchy with self-consistent sources}

In this part, we will construct the super gBK hierarchy with self-consistent sources. At the super-isospectral problem
\eqa  \phi_x=M\phi,\quad \phi_t=N\phi. \label{3.1}\eeqa
Let $\lambda=\lambda_j$, the spectral vector corresponding $\phi$ remember to $\phi_j$, we obtain the the linear system as following
\eqa \left(\begin{array}{c}
\phi_{1j}\\
\phi_{2j}\\
\phi_{3j}\end{array} \right)_x=M_j\left(\begin{array}{c}
\phi_{1j}\\
\phi_{2j}\\
\phi_{3j}\end{array} \right),\quad\left(\begin{array}{c}
\phi_{1j}\\
\phi_{2j}\\
\phi_{3j}\end{array} \right)_t=N_j\left(\begin{array}{c}
\phi_{1j}\\
\phi_{2j}\\
\phi_{3j}\end{array} \right),\label{3.2}\eeqa
where  $M_j=M|_{\lambda=\lambda_j}$, $N_j=N|_{\lambda=\lambda_j}$, $j=1,2...N$. By
\eqa  \frac{\delta \tilde{H}_n}{\delta
u}=\sum^N_{j=1}\frac{\delta \lambda_j}{\delta
u}=\sum^N_{j=1}\left(\begin{array}{c}
Str(\Psi_j\frac{\delta M}{\delta v})\\
Str(\Psi_j\frac{\delta M}{\delta w})\\
Str(\Psi_j\frac{\delta M}{\delta \alpha})\\
Str(\Psi_j\frac{\delta M}{\delta \beta})\end{array}\right)=\left(\begin{array}{c}
<\Phi_1,\Phi_2>\\
2<\Phi_1,\Phi_1>\\
-2<\Phi_2,\Phi_3>\\
2<\Phi_1,\Phi_3>\end{array}\right),\label{3.3}\eeqa
where $\Phi_j=(\phi_{j1} ,\cdots,\phi_{jN})^T$, $j=1,2,3$.
So the super gBK hierarchy with self-consistent sources is proposed
\eqa
\small u_t=\left(\begin{array}{c}
v\\
w\\
\alpha\\
\beta\end{array} \right)_t=J\left(\begin{array}{c}
a_{n}\\
-2b_{n}\\
-2\delta_{n}\\
2\rho_{n}\end{array} \right)+J\left(\begin{array}{c}
<\Phi_1,\Phi_2>\\
2<\Phi_1,\Phi_1>\\
-2<\Phi_2,\Phi_3>\\
2<\Phi_1,\Phi_3>
\end{array}\right),\label{3.4} \eeqa
where
\eqa \left\{ \begin{array}{l}
\phi_{1j,x}=\lambda\phi_{1j}+(w-\frac{1}{2}v)\phi_{2j}+\alpha\phi_{3j},\\
\phi_{2j,x}=2v\phi_{1j}-\lambda\phi_{2j}+\beta\phi_{3j},\\
\phi_{3j,x}=\beta\phi_{1j}-\alpha\phi_{2j}.
\end{array}
\right. \label{3.5} \eeqa
For ~$n=2$, we obtain the super gBK hierarchy with self-consistent sources
\eqa
&&v_{t2}=v_{xx}-2vv_x-4w_x-4\alpha\alpha_{xx}+4\alpha_x\beta+4\alpha\beta_x-2\p\sum^N_{j=1}\phi_{1j}^2,\nn\\
&&w_{t2}=-w_{xx}-2(wv)_x-2v_x+2(2w+2)\alpha_x+2w_x\alpha-(2w\beta+\frac{1}{2}\beta v^2)(1+\alpha)-2\beta\beta_x\nn\\
&&\quad\quad\;\;-\p\sum^N_{j=1}\phi_{1j}\phi_{2j}+\alpha\sum^N_{j=1}\phi_{2j}\phi_{3j}+\beta\sum^N_{j=1}\phi_{1j}\phi_{3j},\nn\\
&&\alpha_{t2}=2\alpha_{xx}-2v\alpha_x+2\beta_x-\frac{1}{2}\alpha v_x-\alpha\sum^N_{j=1}\phi_{1j}^2-\sum^N_{j=1}\phi_{1j}\phi_{3j},
\nn\\
&&\beta_{t2}=-2\beta_{xx}-2v\beta_x+2\beta\alpha\alpha_x+2(2w+2)\alpha_x+2\alpha w_x-\frac{3}{2}\beta v_x+\beta\sum^N_{j=1}\phi_{1j}^2+\sum^N_{j=1}\phi_{2j}\phi_{3j}.\nn\eeqa
\vskip .3in

\section{ Conservation laws for the super gBK hierarchy}

In the following, we will construct conservation laws of the super gBK hierarchy. Introducing two variables
\eqa F=\dfrac{\phi_2}{\phi_1},\quad G=\dfrac{\phi_3}{\phi_1}.\label{4.1}\eeqa
So, we have
\eqa
&& F_x=-2w-2+(2\lambda-v) F+\beta G-F^2-\alpha FG,\label{4.2}\\
&& G_x=\beta-\alpha F+(\lambda-\frac{1}{2}v) G-GF-\alpha G^2.\label{4.3}\eeqa
We expand ~$F,G$ in powers of ~$\lambda^{-1}$
as follows
\eqa F=\sum^\infty_{j=1}f_j\lambda^{-j},\quad G=\sum^\infty_{j=1}g_j\lambda^{-j}.\label{4.4}\eeqa
Substituting Eq.\eqref{4.4} into Eq.\eqref{4.2},\eqref{4.3} and comparing the coefficients of the same power of $\lambda$, we obtain
\eqa &&\lambda^0:\quad f_1=1+w,g_1=-\beta,\nn\\
&&\lambda^{-1}:\; f_2=\frac{1}{2}f_{1x}+\frac{1}{2}vf_{1}-\frac{1}{2}\beta g_{1}=\frac{1}{2}w_x+\frac{1}{2}v(1+w),\nn\\
&& \quad\quad\;\; g_2=g_{1x}+\alpha f_1+\frac{1}{2}vg_1=-\beta_x+\alpha(1+w)-\frac{1}{2}v\beta,\nn\\
&&\lambda^{-2}:\;f_3=\frac{1}{2}f_{2x}+\frac{1}{2}vf_{2}-\frac{1}{2}\beta g_2+\frac{1}{2}f_1^2+\frac{1}{2}\alpha f_1g_1\nn\\&&\quad\quad\quad\;\;
=\frac{1}{4}w_{xx}+\frac{1}{8}w_x^2+\frac{1}{2}vw_x+\frac{1}{4}(1+w)(w_xv+v_x)+\frac{1}{8}v^2(1+w)(3+w)+\frac{1}{2}\beta\beta_x,\nn\\&&\quad\;\;\;\;\;
g_3=g_{2x}+\alpha f_2+\frac{1}{2}vg_2+f_1g_1+\alpha g_1^2\nn\\&&\quad\quad\quad\;\;
=-\beta_{xx}-v\alpha_x+(\alpha_x+\alpha v-\beta)(1+w)+\frac{3}{2}\alpha w_x-\frac{1}{2}v_x\beta-\frac{1}{4}v^2\beta.\nn\eeqa
and the recursion formulas for ~$f_n$ and ~$g_n$ are given
\eqa &&
f_{n+1}=\frac{1}{2}f_{nx}+\frac{1}{2}vf_{n}-\frac{1}{2}\beta g_n+\frac{1}{2}\sum^n_{l=1}f_lf_{n-l}+\frac{1}{2}\alpha\sum^n_{l=1}f_lg_{n-l},\label{4.5}\\
&&g_{n+1}=g_{nx}+\alpha f_n+\frac{1}{2}vg_n+\sum^n_{l=1}f_lg_{n-l}.\label{4.6}\eeqa
Because of linear spectral problems Eq.\eqref{3.1}
\eqa
&&(ln\phi_1)_x=-\lambda+\frac{1}{2}v+F+\alpha G,\label{4.7}\\
&&(ln\phi_1)_t=A+BF+\rho G.
\label{4.8}\eeqa
It is easy to calculate that
\eqa \frac{\partial}{\partial t}((-\lambda+\frac{1}{2}v+F+\alpha G)=\frac{\partial}{\partial x}(A+BF+\rho G).\label{4.9}\eeqa
where
\eqa A=k_0(\lambda^2+w-\alpha\beta),B=-k_0(\lambda+\frac{1}{2}v),\rho=-k_0(\alpha\lambda-\alpha_x+\frac{1}{2}v\alpha).\label{4.10}\eeqa
In order to obtain the conservation laws for super integrable hierarchy, we difine
\eqa \sigma=-\lambda+\frac{1}{2}v+F+\alpha G,\quad \theta=A+BF+\rho G.\nn\eeqa
Then the Eq.\eqref{4.9} can be rewritten as $\sigma_t=\theta_x$, which is just the formal definition of conservation laws. We expand ~$\sigma$ and $\theta$ as series in powers of $\lambda$ with the coefficients, which are called conserved densities and fluxes respectively
~\eqa &&\sigma=-\lambda+\frac{1}{2}v+\sum^\infty_{j=1}\sigma_j\lambda^{-j},\nn\\
&&\theta=k_0\lambda^2-2+\sum^\infty_{j=1}\theta_j\lambda^{-j}.\label{4.11}\eeqa
The first of the conservation of density and flow
~\eqa &&\sigma_1=f_1+\alpha g_1=w+1-\alpha\beta,\nn\\
&&\theta_1=k_0[-f_2-\frac{1}{2}vf_1-\alpha g_2+(\alpha_x-\frac{1}{2}v\alpha)g_1]\nn\\&&\quad\;
=k_0[-\frac{1}{2}w_x-v(1+w)+\alpha\beta_x-\alpha_x\beta+v\alpha\beta]
.\nn\eeqa
With the help of Eq.\eqref{4.9},\eqref{4.10},\eqref{4.11}, the recursion relation for $a_n$ and ~$\theta_n$ are given
\eqa &&\sigma_n=f_n+\alpha g_n,\nn\\
&&\theta_n=k_0[-f_{n+1}-\frac{1}{2}vf_n-\alpha g_{n+1}+(\alpha_x-\frac{1}{2}v\alpha)g_n]
.\nn\eeqa
where  $f_n$ and  $g_n$ can be calculated from Eq.\eqref{4.5} and Eq.\eqref{4.6}.
\vskip .3in

\section{ The novel symmetry constraint }

In order to compute a symmetry constraint, we consider the spectral
problem in Eq.(2.1) and its adjoint spectral problem
 \eqa &&\psi_x=-M^{St}\psi=\left(\begin{array}{ccc}
\lambda-\frac{v}{2}&2w+2&\beta\\
-1&-\lambda+\frac{v}{2}&-\alpha\\
-\alpha&-\beta&0\end{array} \right)\psi,\quad
\psi=\left(\begin{array}{c}
\psi_1\\
\psi_2\\
\psi_3\end{array} \right).\label{5.1}\eeqa
where \lq\lq St \rq\rq means the super transposition. The following result is a general formula for
the variational derivative with respect to the potential $u$ (see \cite{ma1994} for the classical case).

\begin{lem}
(see \cite{he2008,yu2009}) Let $M(u,\lambda)$ be an even matrix of order $m+n$ depending on $u,u_x,u_{xx},\cdots$ and a parameter $\lambda$.
 Suppose that $\phi=(\phi_e,\phi_o)^T$ and $\psi=(\psi_e,\psi_o)^T$ satisfy the spectral problem and the adjoint spectral problem
\beq \phi_x=M(u,\lambda)\phi, \quad \psi_x=-M(u,\lambda)^{St}\psi, \label{5.2} \eeq
where $\phi_e=(\phi_1,\cdots,\phi_m)$ and $\psi_e=(\psi_1,\cdots,\psi_m)$ are even eigenfunctions, and
$\phi_o=(\phi_{m+1},\cdots,\phi_{m+n})$ and $\psi_o=(\psi_{m+1},\cdots,\psi_{m+n})$ are odd eigenfunctions.
Then the variational derivative of the parameter $\lambda$ with respect to the potential $u$ is given by
\beq \frac{\delta\lambda}{\delta u}=\frac{(\psi_e,(-1)^{p(u)}\psi_o)(\frac{\p M}{\p u})\phi}{-\int\psi^T(\frac{\p M}{\p u})\phi dx}, \label{5.3}\eeq
where we denote
\beq p(v)=\left\{\begin{array}{l}
0, \quad $v is an even variable$,\\
1, \quad $v is an odd variable$.
\end{array}\right. \label{5.4}\eeq
\end{lem}

By Lemma 5.1, it's easy to get the variational derivative of the
spectral parameter $\lambda$ with respect to the potential $u$
\beq \frac{\delta\lambda}{\delta u}
 =\left(\begin{array}{c}
\frac{\delta\lambda}{\delta v}\\
\frac{\delta\lambda}{\delta w}\\
\frac{\delta\lambda}{\delta \alpha}\\
\frac{\delta\lambda}{\delta \beta}\end{array}\right)
=\frac{1}{E}\left(\begin{array}{c}
\frac{1}{2}(\psi_1\phi_1-\psi_2\phi_2)\\
-2\psi_2\phi_1\\
\psi_3\phi_2+\psi_1\phi_3\\
\psi_2\phi_3-\psi_3\phi_1
\end{array}
\right), \label{5.5}\eeq
where $E=\int_{-\infty}^{\infty}(\psi_1\phi_1-\psi_2\phi_2)dx$. When zero boundary conditions
$\lim_{|x|\rightarrow\infty}\phi=\lim_{|x|\rightarrow \infty}\psi=0$
are imposed, we can verify a simple characteristic property of the variational derivative
 \beq \mathcal{L}\frac{\delta\lambda}{\delta u}=\lambda\frac{\delta\lambda}{\delta u}, \label{5.6} \eeq
where $\mathcal{L}$ is defined in \eqref{2.6}. Consider the spatial system
 \beq \left\{\begin{array}{l}
\left(\begin{array}{c}
\phi_{1j,x}\\
\phi_{2j,x}\\
\phi_{3j,x}\end{array} \right)=\left(\begin{array}{ccc}
-\lambda+\frac{v}{2}&1&\alpha\\
-2w-2&\lambda-\frac{v}{2}&\beta\\
\beta&-\alpha&0
\end{array} \right)\left(\begin{array}{c}
\phi_{1j}\\
\phi_{2j}\\
\phi_{3j}\end{array} \right),\\
\left(\begin{array}{c}
\psi_{1j,x}\\
\psi_{2j,x}\\
\psi_{3j,x}\end{array} \right)=\left(\begin{array}{ccc}
\lambda-\frac{v}{2}&2w+2&\beta\\
-1&-\lambda+\frac{v}{2}&-\alpha\\
-\alpha&-\beta&0\end{array} \right)\left(\begin{array}{c}
\psi_{1j}\\
\psi_{2j}\\
\psi_{3j}\end{array} \right),\\
\end{array}
\right. \label{5.7}\eeq
and the temporal system
{\beq\left\{\begin{array}{l}
\left(\begin{array}{c}
\phi_{1j}\\
\phi_{2j}\\
\phi_{3j}\end{array}\right)_{t_n}
=\left(\begin{array}{ccc}
\dsum_{i=0}^na_i\lambda_j^{n-i}+b_{n+1} & \dsum_{i=0}^nb_i\lambda_j^{n-i} & \dsum_{i=0}^n\rho_i\lambda_j^{n-i} \\
\dsum_{i=0}^n(-2b_i+2c_i)\lambda_j^{n-i} & -\dsum_{i=0}^na_i\lambda_j^{n-i}-b_{n+1} & \dsum_{i=0}^n\delta_i\lambda_j^{n-i} \\
\dsum_{i=0}^n\delta_i\lambda_j^{n-i} & -\dsum_{i=0}^n\rho_i\lambda_j^{n-i} & 0 \\
\end{array}\right)\left(\begin{array}{c}\phi_{1j}\\\phi_{2j}\\\phi_{3j}\end{array}\right),\\\\
\left(\begin{array}{c}
\psi_{1j}\\
\psi_{2j}\\
\psi_{3j}\end{array}\right)_{tn}
=\left(\begin{array}{ccc}
-\dsum_{i=0}^na_i\lambda_j^{n-i}-b_{n+1} & \dsum_{i=0}^n(2b_i-2c_i)\lambda_j^{n-i} & \dsum_{i=0}^n\delta_i\lambda_j^{n-i} \\
-\dsum_{i=0}^nb_i\lambda_j^{n-i} & \dsum_{i=0}^na_i\lambda_j^{n-i}+b_{n+1} & -\dsum_{i=0}^n\rho_i\lambda_j^{n-i} \\
-\dsum_{i=0}^n\rho_i\lambda_j^{n-i} & -\dsum_{i=0}^n\delta_i\lambda_j^{n-i} & 0 \\
\end{array}\right)\left(\begin{array}{c}\psi_{1j}\\\psi_{2j}\\\psi_{3j}\end{array}\right),
\end{array}\right.\label{5.8}\eeq}
where $\{\lambda_j,j=1,2,\cdots,N\}$ are N distinct eigenparameters, $\{\phi_j\}$ and $\{\psi_j\}$ are corresponding
eigenfunctions and adjoint eigenfunctions, $1\leq j\leq N$.
Now for Eq.\eqref{5.7} and Eq.\eqref{5.8}, we have the following symmetry constraints:
\beq \frac{\delta\tilde{H_k}}{\delta u}=\sum_{j=1}^N \gamma_j\frac{\delta\lambda_j}{\delta u},\quad k\geq0, \label{5.9}\eeq
where letting $\gamma_j=E_j=\int_{-\infty}^{\infty}(\phi_{1j}\psi_{1j}-\phi_{2j}\psi_{2j})dx$.
If letting initial value $a_0=k_0=1$, and seek $k=1$, we have the following novel symmetry constraint
\beq \left\{\begin{array}{l}
w-\alpha\beta=\frac{1}{2}(\langle\Psi_1,\Phi_1\rangle-\langle\Psi_2,\Phi_2\rangle), \\
v=-2\langle\Psi_2,\Phi_1\rangle, \\
-2\beta_x-v\beta=\langle\Psi_3,\Phi_2\rangle+\langle\Psi_1,\Phi_3\rangle,\\
-2\alpha_x+v\alpha=\langle\Psi_2,\Phi_3\rangle-\langle\Psi_3,\Phi_1\rangle.
\end{array}
\right. \label{5.10}\eeq
where we use the following notation
$$\Psi_i=(\psi_{i1},\cdots,\psi_{iN})^T,\quad \Phi_i=(\phi_{i1},\cdots,\phi_{iN})^T,\quad i=1,2,3,$$
and $\langle\cdot,\cdot\rangle$ denotes the standard inner product of the Euclidian space $R^N$.
We find that the even potentials $v$ and $w$ can be explicitly
expressed by eigenfunctions, but the odd potentials $\alpha$ and $\beta$ cannot. So the symmetry
constraint \eqref{5.10} is called a novel constraint.
\vskip .3in

\section{ The Binary nonlinearization}

In this part, we introduce the following new independent odd variables as a result of the odd potentials
$\alpha$ and $\beta$ cannot be explicitly expressed by eigenfunctions
\beq \phi_{N+1}=\alpha,\quad \psi_{N+1}=2\beta. \label{6.1}\eeq
Considering the new variables of Eq.\eqref{4.1} and substituting Eq.\eqref{5.10} into Eq.\eqref{5.7}, we have the following finite-dimensional super system
{\beq\left\{\begin{array}{l}
\phi_{1j,x}=(-\lambda_j-\langle\Psi_2,\Phi_1\rangle)\phi_{1j}+\phi_{2j}+\phi_{N+1}\phi_{3j}, \\
\phi_{2j,x}=-(\langle\Psi_1,\Phi_1\rangle-\langle\Psi_2,\Phi_2\rangle+\phi_{N+1}\psi_{N+1}+2)\phi_{1j}
+(\lambda_j+\langle\Psi_2,\Phi_1\rangle)\phi_{2j}+\frac{1}{2}\psi_{N+1}\phi_{3j},\\
\phi_{3j,x}=\frac{1}{2}\psi_{N+1}\phi_{1j}-\phi_{N+1}\phi_{2j},\\
\phi_{N+1,x}=-\frac{1}{2}(\langle\Psi_2,\Phi_3\rangle-\langle\Psi_3,\Phi_1\rangle)-\langle\Psi_2,\Phi_1\rangle\phi_{N+1},\\
\psi_{1j,x}=(\lambda_j+\langle\Psi_2,\Phi_1\rangle)\psi_{1j}+(\langle\Psi_1,\Phi_1\rangle-\langle\Psi_2,\Phi_2\rangle+\phi_{N+1}\psi_{N+1}+2)\psi_{2j}
+\frac{1}{2}\psi_{N+1}\psi_{3j},\\
\psi_{2j,x}=-\psi_{1j}+(-\lambda_j-\langle\Psi_2,\Phi_1\rangle)\psi_{2j}-\phi_{N+1}\psi_{3j}, \\
\psi_{3j,x}=-\phi_{N+1}\psi_{1j}-\frac{1}{2}\psi_{N+1}\psi_{2j},\\
\psi_{N+1,x}=-(\langle\Psi_3,\Phi_2\rangle+\langle\Psi_1,\Phi_3\rangle)+\langle\Psi_2,\Phi_1\rangle\psi_{N+1},\\
\end{array}
\right. \label{6.2}\eeq }
where $\Lambda=\diag(\lambda_1,\lambda_2,\cdots,\lambda_N)$, obviously, the system \eqref{6.2} can be written to the following Hamiltonian form
\beq \left\{\begin{array}{l}
\Phi_{1,x}=\frac{\p H_1}{\p \Psi_{1}},\Phi_{2,x}=\frac{\p H_1}{\p\Psi_{2}},
\Phi_{3,x}=\frac{\p H_1}{\p \Psi_{3}},\phi_{N+1,x}=\frac{\p H_1}{\p \psi_{N+1}},\\
\Psi_{1,x}=-\frac{\p H_1}{\p \Phi_{1}},\Psi_{2,x}=-\frac{\p H_1}{\p \Phi_{2}},
\Psi_{3,x}=\frac{\p H_1}{\p \Phi_{3}},\psi_{N+1,x}=\frac{\p H_1}{\p \phi_{N+1}},\\
\end{array} \right. \label{6.3}\eeq
where
\eqa&&
H_1=-\langle\Lambda\Psi_1,\Phi_1\rangle+\langle\Lambda\Psi_2,\Phi_2\rangle-2\langle\Psi_2,\Phi_1\rangle+\langle\Psi_1,\Phi_2\rangle\nn\\&& \quad\quad\;\;
-\langle\Psi_2,\Phi_1\rangle(\langle\Psi_1,\Phi_1\rangle-\langle\Psi_2,\Phi_2\rangle)
-\phi_{N+1}\psi_{N+1}\langle\Psi_2,\Phi_1\rangle\nn\\&& \quad\quad\;\;
+\phi_{N+1}(\langle\Psi_3,\Phi_2\rangle+\langle\Psi_1,\Phi_3\rangle)
+\frac{1}{2}\psi_{N+1}(\langle\Psi_2,\Phi_3\rangle-\langle\Psi_3,\Phi_1\rangle). \nn\eeqa
As for the $t_2$-part, substituting the symmetry constraint \eqref{5.10} into system \eqref{5.8}, we obtain the following finite-dimensional system
{\beq\left\{\begin{array}{l}
\phi_{1j,t_2}=(\lambda^2+\frac{1}{4}(\tilde{v}_{x}-\tilde{v}^2)-\tilde{\alpha}\tilde{\alpha}_x)\phi_{1j}
+(-\lambda-\frac{1}{2}\tilde{v})\phi_{2j}+(-\tilde{\alpha}\lambda+\tilde{\alpha}_x-\frac{1}{2}\tilde{v}\tilde{\alpha})\phi_{3j}, \\
\phi_{2j,t_2}=(2(1+\tilde{w})\lambda+\tilde{w}_x+\tilde{v}(1+\tilde{w}))\phi_{1j}
+(-\lambda^2-\frac{1}{4}(\tilde{v}_{x}-\tilde{v}^2)+\tilde{\alpha}\tilde{\alpha}_x)\phi_{2j}\\ \quad\quad\quad\;\;
+(-\tilde{\beta}\lambda-\tilde{\beta}_x-\frac{1}{2}\tilde{v}\tilde{\beta})\phi_{3j},\\
\phi_{3j,t_2}=(-\tilde{\beta}\lambda-\tilde{\beta}_x-\frac{1}{2}\tilde{v}\tilde{\beta})\phi_{1j}
+(\tilde{\alpha}\lambda-\tilde{\alpha}_x+\frac{1}{2}\tilde{v}\tilde{\alpha})\phi_{2j},\\
\psi_{1j,t_2}=(-\lambda^2-\frac{1}{4}(\tilde{v}_{x}-\tilde{v}^2)+\tilde{\alpha}\tilde{\alpha}_x)\psi_{1j}
+(-2(1+\tilde{w})\lambda-\tilde{w}_x-\tilde{v}(1+\tilde{w}))\psi_{2j}\\ \quad\quad\quad\;\;
+(-\tilde{\beta}\lambda-\tilde{\beta}_x-\frac{1}{2}\tilde{v}\tilde{\beta})\psi_{3j},\\
\psi_{2j,t_2}=(\lambda+\frac{1}{2}\tilde{v})\psi_{1j}+(\lambda^2+\frac{1}{4}(\tilde{v}_{x}-\tilde{v}^2)-\tilde{\alpha}\tilde{\alpha}_x)\psi_{2j}
+(\tilde{\alpha}\lambda-\tilde{\alpha}_x+\frac{1}{2}\tilde{v}\tilde{\alpha})\psi_{3j}, \\
\psi_{3j,t_2}=(\tilde{\alpha}\lambda-\tilde{\alpha}_x-\frac{1}{2}\tilde{v}\tilde{\alpha})\psi_{1j}
+(\tilde{\beta}\lambda+\tilde{\beta}_x+\frac{1}{2}\tilde{v}\tilde{\beta})\psi_{2j},\\
\end{array}\right. \label{6.4}\eeq}
where $\tilde{v},\tilde{w},\tilde{\alpha},\tilde{\beta}$ denote the functions $v,w,\alpha,\beta$
defined by the symmetry constraint given in\eqref{5.10} and
$\tilde{v}_{x},\tilde{w}_{x},\tilde{\alpha}_x,\tilde{\beta}_x$ are given by the following identities:
\eqa &&\tilde{v}_{x}=4\langle\Lambda\Psi_2,\Phi_1\rangle+4\langle\Psi_2,\Phi_1\rangle^2+
2(\langle\Psi_1,\Phi_1\rangle-\langle\Psi_2,\Phi_2\rangle)-2\phi_{N+1}(\langle\Psi_2,\Phi_3\rangle-\langle\Psi_3,\Phi_1\rangle), \nn\\&&
\tilde{w}_{x}=\langle\Psi_1,\Phi_2\rangle+(\langle\Psi_1,\Phi_1\rangle-\langle\Psi_2,\Phi_2\rangle+\phi_{N+1}\psi_{N+1}+2)\langle\Psi_2,\Phi_1\rangle
, \nn\\&&
\tilde{\alpha}_x=-\frac{1}{2}(\langle\Psi_2,\Phi_3\rangle-\langle\Psi_3,\Phi_1\rangle)-\langle\Psi_2,\Phi_1\rangle\phi_{N+1}, \nn\\&&
\tilde{\beta}_x=-\frac{1}{2}(\langle\Psi_3,\Phi_2\rangle+\langle\Psi_1,\Phi_3\rangle)+\frac{1}{2}\langle\Psi_2,\Phi_1\rangle\psi_{N+1}
.\nn\eeqa
Considering the new variables in Eq.\eqref{6.1} and the symmetry constraint given in Eq.\eqref{5.10},
the above finite-dimensional super system given in Eq.\eqref{6.4} becomes the following finite-dimensional system
{\beq\left\{\begin{array}{l}
 \phi_{1j,t_2}=[\lambda^2+\langle\Lambda\Psi_2,\Phi_1\rangle+\frac{1}{2}(\langle\Psi_1,\Phi_1\rangle-\langle\Psi_2,\Phi_2\rangle)]\phi_{1j}\\ \quad\quad\quad\;\;
+(-\lambda+\langle\Psi_2,\Phi_1\rangle)\phi_{2j}+[-\phi_{N+1}\lambda-\frac{1}{2}(\langle\Psi_2,\Phi_3\rangle-\langle\Psi_3,\Phi_1\rangle)]\phi_{3j}, \\
  \phi_{2j,t_2}=[(2+(\langle\Psi_1,\Phi_1\rangle-\langle\Psi_2,\Phi_2\rangle)+\phi_{N+1}\psi_{N+1})\lambda+\langle\Psi_1,\Phi_2\rangle]\phi_{1j}\\
\quad\quad\quad\;\;+[-\lambda^2-\langle\Lambda\Psi_2,\Phi_1\rangle-\frac{1}{2}(\langle\Psi_1,\Phi_1\rangle-\langle\Psi_2,\Phi_2\rangle)]\phi_{2j}\\ \quad\quad\quad\;\;+[-\frac{1}{2}\psi_{N+1}\lambda+\frac{1}{2}(\langle\Psi_3,\Phi_2\rangle+\langle\Psi_1,\Phi_3\rangle)]\phi_{3j},\\
  \phi_{3j,t_2}=[-\frac{1}{2}\psi_{N+1}\lambda+\frac{1}{2}(\langle\Psi_3,\Phi_2\rangle+\langle\Psi_1,\Phi_3\rangle)]\phi_{1j}
+[\phi_{N+1}\lambda+\frac{1}{2}(\langle\Psi_2,\Phi_3\rangle-\langle\Psi_3,\Phi_1\rangle)]\phi_{2j},\\
  \phi_{N+1,t_2}=-\frac{1}{2}(\langle\Lambda\Psi_2,\Phi_3\rangle-\langle\Lambda\Psi_3,\Phi_1\rangle)-\phi_{N+1}\langle\Lambda\Psi_2,\Phi_1\rangle,\\
  \psi_{1j,t_2}=[-\lambda^2-\langle\Lambda\Psi_2,\Phi_1\rangle-\frac{1}{2}(\langle\Psi_1,\Phi_1\rangle-\langle\Psi_2,\Phi_2\rangle)]\psi_{1j}\\ \quad\quad\quad\;\;
+[(-2-(\langle\Psi_1,\Phi_1\rangle+\langle\Psi_2,\Phi_2\rangle)+\phi_{N+1}\psi_{N+1})\lambda-\langle\Psi_1,\Phi_2\rangle]\psi_{2j}\\ \quad\quad\quad\;\;
+[-\frac{1}{2}\psi_{N+1}\lambda+\frac{1}{2}(\langle\Psi_3,\Phi_2\rangle+\langle\Psi_1,\Phi_3\rangle)]\psi_{3j},\\
  \psi_{2j,t_2}=(\lambda-\langle\Psi_2,\Phi_1\rangle)\psi_{1j}+[\lambda^2+\langle\Lambda\Psi_2,\Phi_1\rangle+\frac{1}{2}(\langle\Psi_1,\Phi_1\rangle-\langle\Psi_2,\Phi_2\rangle)]\psi_{2j}\\ \quad\quad\quad\;\;
+[\phi_{N+1}\lambda+\frac{1}{2}(\langle\Psi_2,\Phi_3\rangle-\langle\Psi_3,\Phi_1\rangle)]\psi_{3j}, \\
  \psi_{3j,t_2}=[\phi_{N+1}\lambda+\frac{1}{2}(\langle\Psi_2,\Phi_3\rangle-\langle\Psi_3,\Phi_1\rangle)]\psi_{1j}
+[\frac{1}{2}\psi_{N+1}\lambda-\frac{1}{2}(\langle\Psi_3,\Phi_2\rangle+\langle\Psi_1,\Phi_3\rangle)]\psi_{2j},\\
\psi_{N+1,t_2}=-(\langle\Lambda\Psi_3,\Phi_2\rangle+\langle\Lambda\Psi_1,\Phi_3\rangle)+\psi_{N+1}\langle\Lambda\Psi_2,\Phi_1\rangle.\\
\end{array}\right. \label{6.5}\eeq }
By a direct but tedious calculation, the finite-dimensional system \eqref{6.5} become to the following Hamiltonian form
\beq \left\{\begin{array}{l}
\Phi_{1,t_2}=\frac{\p H_2}{\p \Psi_{1}},\Phi_{2,t_2}=\frac{\p H_2}{\p\Psi_{2}},
\Phi_{3,t_2}=\frac{\p H_2}{\p \Psi_{3}},\phi_{N+1,t_2}=\frac{\p H_2}{\p \psi_{N+1}},\\
\Psi_{1,t_2}=-\frac{\p H_2}{\p \Phi_{1}},\Psi_{2,t_2}=-\frac{\p H_2}{\p \Phi_{2}},
\Psi_{3,t_2}=\frac{\p H_2}{\p \Phi_{3}},\psi_{N+1,t_2}=\frac{\p H_2}{\p \phi_{N+1}},
\end{array} \right. \label{6.6}\eeq
where the Hamilton function is
\eqa
&&H_2=\langle\Lambda^2\Psi_1,\Phi_1\rangle-\langle\Lambda^2\Psi_2,\Phi_2\rangle+\langle\Lambda\Psi_2,\Phi_1\rangle(\langle\Psi_1,\Phi_1\rangle-\langle\Psi_2,\Phi_2\rangle)
+2\langle\Lambda\Psi_2,\Phi_1\rangle
\nn\\&&\quad\quad\;-\langle\Lambda\Psi_1,\Phi_2\rangle+\langle\Psi_2,\Phi_1\rangle\langle\Psi_1,\Phi_2\rangle-\phi_{N+1}(\langle\Lambda\Psi_1,\Phi_3\rangle+\langle\Lambda\Psi_3,\Phi_2\rangle)-\frac{1}{2}(\langle\Psi_2,\Phi_3\rangle
\nn\\&&\quad\quad\;-\langle\Psi_3,\Phi_1\rangle)(\langle\Psi_3,\Phi_2\rangle+\langle\Psi_1,\Phi_3\rangle)+\frac{1}{2}(\langle\Lambda\Psi_2,\Phi_3\rangle-\langle\Lambda\Psi_3,\Phi_1\rangle)\psi_{N+1}
\nn\\&&\quad\quad\;+\phi_{N+1}\psi_{N+1}\langle\Lambda\Psi_2,\Phi_1\rangle+\frac{1}{4}(\langle\Psi_1,\Phi_1\rangle-\langle\Psi_2,\Phi_2\rangle)^{2}.
\nn\eeqa
In the following, we prove that for any $n\geq2$, the super system given in Eq.\eqref{5.8} can be nonlinearized
and furthermore, the obtained nonlinearized system is a finite-dimensional super-Hamiltonian system. Therefore,
making use of Eq.\eqref{5.6} and the recursion relation \eqref{2.5}, yields
\beq \left\{
\begin{array}{l}
\tilde{a}_{m+1}=\frac{1}{2}(\langle\Lambda^{m-1}\Psi_1,\Phi_1\rangle-\langle\Lambda^{m-1}\Psi_2,\Phi_2\rangle),\quad m\geq1,\\
\tilde{b}_{m+1}=\langle\Lambda^{m-1}\Psi_2,\Phi_1\rangle, \quad m\geq1,\\
\tilde{c}_{m+1}=\langle\Lambda^{m-1}\Psi_2,\Phi_1\rangle+\frac{1}{2}\langle\Lambda^{m-1}\Psi_1,\Phi_2\rangle,\quad m\geq1,\\
\tilde{\rho}_{m+1}=-\frac{1}{2}(\langle\Lambda^{m-1}\Psi_2,\Phi_3\rangle-\langle\Lambda^{m-1}\Psi_3,\Phi_1\rangle),\quad m\geq1,\\
\tilde{\delta}_{m+1}=\frac{1}{2}(\langle\Lambda^{m-1}\Psi_3,\Phi_2\rangle+\langle\Lambda^{m-1}\Psi_1,\Phi_3\rangle),\quad m\geq1,\\
\end{array}
\right. \label{6.7}\eeq
where $\Lambda=$diag$(\lambda_1,\lambda_2,\cdots,\lambda_N)$. Substituting Eq.\eqref{6.7} into Eq.\eqref{5.8}, we have
{\beq\left\{\begin{array}{l}
\phi_{1j,t_n}=(\dsum_{i=0}^n\tilde{a}_i\lambda_j^{n-i}+\tilde{b}_{n+1})\phi_{1j}
+\dsum_{i=0}^n\tilde{b}_i\lambda_j^{n-i}\phi_{2j}+\dsum_{i=0}^n\tilde{\rho}_i\lambda_j^{n-i}\phi_{3j}, \\
\phi_{2j,t_n}=(\dsum_{i=0}^n(-2\tilde{b}_i+2\tilde{c}_i)\lambda_j^{n-i})\phi_{1j}
+(-\dsum_{i=0}^n\tilde{a}_i\lambda_j^{n-i}-\tilde{b}_{n+1})\phi_{2j}
+\dsum_{i=0}^n\tilde{\delta}_i\lambda_j^{n-i}\phi_{3j},\\
\phi_{3j,t_n}=\dsum_{i=0}^n\tilde{\delta}_i\lambda_j^{n-i}\phi_{1j}-\dsum_{i=0}^n\tilde{\rho}_i\lambda_j^{n-i}\phi_{2j},\\
\phi_{N+1,t_n}=\frac{1}{2}(\langle\Lambda^{n-1}\Psi_2,\Phi_3\rangle-\langle\Lambda^{n-1}\Psi_3,\Phi_1\rangle)
-\phi_{N+1}\langle\Lambda^{n-1}\Psi_2,\Phi_1\rangle,\\
\psi_{1j,t_n}=(-\dsum_{i=0}^n\tilde{a}_i\lambda_j^{n-i}-\tilde{b}_{n+1})\psi_{1j}
+(\dsum_{i=0}^n(2\tilde{b}_i-2\tilde{c}_i)\lambda_j^{n-i})\psi_{2j}
+\dsum_{i=0}^n\tilde{\delta}_i\lambda_j^{n-i}\psi_{3j},\\
\psi_{2j,t_n}=-\dsum_{i=0}^n\tilde{b}_i\lambda_j^{n-i}\psi_{1j}+(\dsum_{i=0}^n\tilde{a}_i\lambda_j^{n-i}+\tilde{b}_{n+1})\psi_{2j}
-\dsum_{i=0}^n\tilde{\rho}_i\lambda_j^{n-i}\psi_{3j}, \\
\psi_{3j,t_n}=-\dsum_{i=0}^n\tilde{\rho}_i\lambda_j^{n-i}\psi_{1j}-\dsum_{i=0}^n\tilde{\delta}_i\lambda_j^{n-i}\psi_{2j},\\
\psi_{N+1,t_n}=-(\langle\Lambda^{n-1}\Psi_3,\Phi_2\rangle+\langle\Lambda^{n-1}\Psi_1,\Phi_3\rangle)+\psi_{N+1}\langle\Lambda^{n-1}\Psi_2,\Phi_1\rangle.\\
\end{array}\right. \label{6.8}\eeq }
Next, we show that the nonlinearized super system given in Eq.\eqref{6.8} is a finite-dimensional super-Hamiltonian
system. Under the constraint \eqref{5.10}, the identity $(\tilde{N})_x=[\tilde{M},\tilde{N}]$ and
$(\tilde{N}^2)_x=[\tilde{M},\tilde{N}^2]$ are still satisfied, we have
\beq F_x=(\frac{1}{2}Str\tilde{N}^2)_x=\frac{d}{dx}(\tilde{a}^2-2\tilde{b}^2+2\tilde{b}\tilde{c}+2\tilde{\rho}\tilde{\delta})=0. \label{6.9}\eeq
The identity indicates that $F$ is a generating function of integrals of motion for the nonlinearized spatial systems
\eqref{6.2}. Let $F = \sum_{n\geq0} F_n\lambda^{-n}$, and we obtain the following formulas:
\beq F_m=\sum_{i=0}^m(\tilde{a}_i\tilde{a}_{m-i}-2\tilde{b}_i\tilde{b}_{m-i}+2\tilde{b}_i\tilde{c}_{m-i}+2\tilde{\rho}_i\tilde{\delta}_{m-i}). \label{6.10}\eeq
Assume that $\tilde{a}_0=1,\tilde{b}_0=\tilde{c}_0=\tilde{\rho}_0=\tilde{\delta}_0=0$ and from the above relations, we find
\eqa
&&F_0=1,\quad F_1=2\tilde{a}_1=0,\quad F_2=-2,\nn\\
&&F_3=\langle\Lambda\Psi_1,\Phi_1\rangle-\langle\Lambda\Psi_2,\Phi_2\rangle+4\langle\Lambda^{2}\Psi_2,\Phi_1\rangle-2\langle\Psi_2,\Phi_1\rangle-\langle\Psi_1,\Phi_2\rangle\nn\\
&&\quad\quad\;\;+\langle\Psi_2,\Phi_1\rangle(\langle\Psi_1,\Phi_1\rangle-\langle\Psi_2,\Phi_2\rangle+\phi_{N+1}\psi_{N+1})-\phi_{N+1}(\langle\Psi_1,\Phi_3\rangle\nn\\
&&\quad\quad\;\;+\langle\Psi_3,\Phi_2\rangle)+\frac{1}{2}(\langle\Psi_1,\Phi_3\rangle-\langle\Psi_3,\Phi_2\rangle)\psi_{N+1}\nn\\
&&F_4=\langle\Lambda^2\Psi_1,\Phi_1\rangle-\langle\Lambda^2\Psi_2,\Phi_2\rangle+\langle\Lambda\Psi_2,\Phi_1\rangle(\langle\Psi_1,\Phi_1\rangle-\langle\Psi_2,\Phi_2\rangle)
+2\langle\Lambda\Psi_2,\Phi_1\rangle
\nn\\&&\quad\quad\;-\langle\Lambda\Psi_1,\Phi_2\rangle+\langle\Psi_2,\Phi_1\rangle\langle\Psi_1,\Phi_2\rangle-\phi_{N+1}(\langle\Lambda\Psi_1,\Phi_3\rangle+\langle\Lambda\Psi_3,\Phi_2\rangle)-\frac{1}{2}(\langle\Psi_2,\Phi_3\rangle
\nn\\&&\quad\quad\;-\langle\Psi_3,\Phi_1\rangle)(\langle\Psi_3,\Phi_2\rangle+\langle\Psi_1,\Phi_3\rangle)+\frac{1}{2}(\langle\Lambda\Psi_2,\Phi_3\rangle-\langle\Lambda\Psi_3,\Phi_1\rangle)\psi_{N+1}
\nn\\&&\quad\quad\;+\phi_{N+1}\psi_{N+1}\langle\Lambda\Psi_2,\Phi_1\rangle+\frac{1}{4}(\langle\Psi_1,\Phi_1\rangle-\langle\Psi_2,\Phi_2\rangle)^{2}.
\nn\\
&&F_m=\sum_{i=2}^{m-2}(\tilde{a}_i\tilde{a}_{m-i}-2\tilde{b}_i\tilde{b}_{m-i}+2\tilde{b}_i\tilde{c}_{m-i}+2\tilde{\rho}_i\tilde{\delta}_{m-i})+2\tilde{a}_0\tilde{a}_m-4\tilde{b}_1\tilde{b}_{m-1}+2\tilde{b}_1\tilde{c}_{m-1}\nn\\
&&\quad\quad\;\;+2\tilde{b}_{m-1}\tilde{c}_{1}+2\tilde{\rho}_1\tilde{\delta}_{m-1}+2\tilde{\rho}_{m-1}\tilde{\delta}_{1}\nn\\
&&\quad\;\;=\langle\Lambda^{m-2}\Psi_1,\Phi_1\rangle-\langle\Lambda^{m-2}\Psi_2,\Phi_2\rangle+2\langle\Lambda^{m-3}\Psi_2,\Phi_1\rangle-\langle\Lambda^{m-3}\Psi_1,\Phi_2\rangle\nn\\
&&\quad\quad\;\;+\langle\Lambda^{m-3}\Psi_2,\Phi_1\rangle(\langle\Psi_1,\Phi_1\rangle-\langle\Psi_2,\Phi_2\rangle+\phi_{N+1}\psi_{N+1})-\phi_{N+1}(\langle\Lambda^{m-3}\Psi_3,\Phi_2\rangle\nn\\
&&\quad\quad\;\;+\langle\Lambda^{m-3}\Psi_1,\Phi_3\rangle)+\frac{1}{2}\psi_{N+1}(\langle\Lambda^{m-3}\Psi_2,\Phi_3\rangle-\langle\Lambda^{m-3}\Psi_3,\Phi_1\rangle)\nn\\
&&\quad\quad\;\;+\sum_{i=2}^{m-2}[\frac{1}{4}(\langle\Lambda^{i-2}\Psi_1,\Phi_1\rangle-\langle\Lambda^{i-2}\Psi_2,\Phi_2\rangle)(\langle\Lambda^{m-i-2}\Psi_1,\Phi_1\rangle-\langle\Lambda^{m-i-2}\Psi_2,\Phi_2\rangle)\nn\\
&&\quad\quad\;\;+\langle\Lambda^{i-2}\Psi_2,\Phi_1\rangle\langle\Lambda^{m-i-2}\Psi_1,\Phi_2\rangle-\langle\Lambda^{i-2}\Psi_2,\Phi_3\rangle+\langle\Lambda^{i-2}\Psi_3,\Phi_1\rangle].
\nn \eeqa
The constrained temporal part of the super gBK equation hierarchy \eqref{5.8} can be rewritten to the following Hamilton form
\beq \left\{
\begin{array}{l}
\Phi_{1,t_n}=\frac{\p F_{n+2}}{\p\Psi_1},\Phi_{2,t_n}=\frac{\p F_{n+2}}{\p\Psi_2},
\Phi_{3,t_n}=\frac{\p F_{n+2}}{\p\Psi_3},\phi_{N+1,t_n}=\frac{\p F_{n+2}}{\p\psi_{N+1}},\\
\Psi_{1,t_n}=-\frac{\p F_{n+2}}{\p\Phi_1},\Psi_{2,t_n}=-\frac{\p F_{n+2}}{\p\Phi_2},
\Psi_{3,t_n}=\frac{\p F_{n+2}}{\p\Phi_3},\psi_{N+1,t_n}=\frac{\p F_{n+2}}{\p\phi_{N+1}},\\
\end{array}
\right. \label{6.11}\eeq
i.e., for any $n$, the constrained temporal part of the super gBK equation hierarchy \eqref{5.8} is the
finite-dimensional super Hamilton hierarchy. As an calculation example, we have the following equation:
 \eqa&&\Phi_{1,t_n}=(\dsum_{i=1}^{n}\tilde{a}_i\Lambda^{n-i}+\tilde{b}_{n+1} ) \Phi_{1}+\dsum_{i=1}^n\tilde{b}_i\Lambda^{n-i}\Phi_{2}+\dsum_{i=1}^n\tilde{\rho}_i\Lambda^{n-i}\Phi_{3}
\nn\\
&&\quad\quad\;=\Lambda^n\Phi_{1}+\langle\Lambda^{n-1}\Psi_2,\Phi_1\rangle\Phi_{1}-\Lambda^{n-1}\Phi_{2}-\phi_{N+1}\Lambda^{n-1}\Phi_{3}\nn\\
&&\quad\quad\quad\;-\frac{1}{2}\sum_{i=2}^{n}(\langle\Lambda^{i-2}\Psi_1,\Phi_1\rangle-\langle\Lambda^{i-2}\Psi_2,\Phi_2\rangle)\Lambda^{n-i}\Phi_1\nn\\
&&\quad\quad\quad\;+\sum_{i=2}^{n}\langle\Lambda^{i-2}\Psi_2,\Phi_1\rangle\Lambda^{n-i}\Phi_2-\frac{1}{2}\sum_{i=2}^{n}(\langle\Lambda^{i-2}\Psi_2,\Phi_3\rangle-\langle\Lambda^{i-2}\Psi_3,\Phi_1\rangle)\Lambda^{n-i}\Phi_3\nn\\
&&\quad\quad\quad\;
+(\langle\Lambda^{i-1}\Psi_2,\Phi_3\rangle-\langle\Lambda^{i-1}\Psi_3,\Phi_1\rangle)\Lambda^{n-i}\Phi_3=\frac{\p F_{n+2}}{\p\Psi_1}. \nn\eeqa
It is not difficult to see that $F_n (n \geq 0)$ are also integrals of motion for the temporal system \eqref{5.8}, i.e.,
\beq\{F_{m},F_{n+2}\}=\frac{\p}{\p t_n}F_{m}=0.\quad m,n\geq0, \label{6.12}\eeq
where the Poisson bracket is defined by
 \eqa&&\{f,g\}=\sum_{i=1}^3\sum_{i=1}^N (\frac{\p f}{\p\phi_{ij}}\frac{\p g}{\p\psi_{ij}}-(-1)^{p(\phi_{ij})p(\psi_{ij})}
\frac{\p f}{\p\psi_{ij}}\frac{\p g}{\p \phi_{ij}})\nn\\
&&\quad\quad\quad
+(\frac{\p f}{\p\phi_{N+1}}\frac{\p g}{\p \psi_{N+1}}
+\frac{\p f}{\p\psi_{N+1}}\frac{\p g}{\p \phi_{N+1}}).\label{6.13}\eeqa
It is natural for us to set
\beq f_k =\phi_{1k}\psi_{1k}+\phi_{2k}\psi_{2k}+\phi_{3k}\psi_{3k}, \quad 1 \leq k \leq N ,\label{6.14}\eeq
and verify they are also integrals of motion of the constrained systems \eqref{6.2} and \eqref{5.8}.
In the same way with \cite{he2008}, we can prove the independence of
$\{f_k\}_{k=1}^N,\{F_k\}_{k=2}^{2N+3}$. So the following theorem hold ture.

\begin{thm}Both the spatial systems \eqref{5.7} and the temporal systems
\eqref{5.8}  under the symmetry constraint \eqref{5.10} become completely
finite-dimensional integrable Hamiltonian systems in the Liouville sense.
\end{thm}
\vskip .3in

\section{Conclusions and discussions}

Starting from Lie super algebras, we may get super equation hierarchy.
With the help of variational identity, the Hamiltonian structure can also be presented.
Based on Lie super algebra, the self-consistent sources of super gBK
hierarchy can be obtained. It enriched the content of self-consistent sources of super
soliton hierarchy. In addition, we also get the conservation laws of the super gBK hierarchy. Finally, we have applied the binary nonlinearization method to the
super gBK hierarchy by the symmetry constraint \eqref{5.10}.
It provides a new and systematic way to construct a finite-dimensional super Hamiltonian
system. The methods in this study can be applied to other super soliton hierarchy to get
more super hierarchies with self-consistent. And under constraints of this form, the
nonlinearization of the other super soliton hierarchy will be studied in our future work.
\\
\begin{rem} With the development of soliton theory, super integrable systems associated with Lie
super algebra have been paid growing attention, many classical integrable equations
have been extended to be the super completely integrable equations. However, according to
the eassies now available, scholars studyed either the self-consistent sources and the conservation laws
or the symmetry constraint and the binary nonlinearization of super integrable systems .
So far, only author aired one of his eassies \cite{Hu2017} in which he studied both theories mentioned above.
And in this paper, there is a difference is that the super generalized Broer-Kaup
equation hierarchy have a novel symmetry constraint compared with paper \cite{Hu2017}.
\end{rem}

\end{document}